\documentclass[a4paper,12pt]{article}
\usepackage{graphicx}
\title{A Comparison Study of Two Methods for Elliptic Boundary Value Problems}
\author{\bf Jian Du$^1$, Shuqiang Wang$^1$, James Glimm$^{1,2}$, Roman Samulyak$^2$
\\ \it
$^1$Department of Applied Mathematics and Statistics,
\\ \it
SUNY at Stony Brook, Stony Brook, NY 11794, USA
\\ \it
$^2$Computational Science Center,
\\ \it
Brookhaven National Laboratory, Upton, NY 11973
}
\begin{document}

\maketitle

\begin{abstract}
    In this paper, we perform a comparison study of two methods (the embedded boundary method and several versions of the mixed finite element method) to solve an elliptic boundary value problem.
\end{abstract}

\section{Introduction}
    The purpose of this paper is to present a comparative study of two popular methods for the solution of elliptic boundary value problem: the embedded boundary method (EBM) and the mixed finite element methods (MFEM). The methods are quite different in their performance characteristics and the mixed finite element methods could use different basis functions. 

    To present our main results in an easily accessible manner, we arrange the results in a table of solution time for comparable accuracy. We find that the EBM is better than lower or the same order accurate MFEM, but perhaps not as good as the higher order accurate MFEM we test here. 

    We observe that no single study of comparison can be definitive, as comparison results may be dependent on the problem chosen, the accuracy desired and comparison method selected. To begin, we distinguish between two not so different kinds of elliptic problems: the elliptic boundary value problems and the elliptic interface problem. For the elliptic boundary value problem, the computational domain exists only on one side of the boundary, for example, interior/exterior boundary value problem. For the elliptic interface problem, there is some internal boundary called an interface across which the solutions on the two sides satisfy some jump conditions. 

    There are many methods for solving the elliptic boundary value/interface problems. Several popular methods have been developed on cartesian meshes for the boundary value/interface problems: the immersed boundary method (IBM) by Peskin \cite{Peskin2002}, the immersed interface method by LeVeque and Li \cite{LeVequeLi1994}, the ghost fluid methods (GFM) by Liu, etc. \cite{LiuFedkiwKang2000}, the embedded boundary method by Johansen and Colella \cite{JohansenColella98}, integral equation method by Mayo \cite{Mayo1984}, Mckenney, Greengard and Mayo \cite{MckenneyGreenGardMayo1995}. The advantage of these methods is that they are defined on a cartesian mesh. Therefore no need to generate a mesh. For the cells away from the boundary/interface, they just use a central finite difference method which is simple and second order accurate. For the cells near or crossing the boundary/interface, special treatment is needed. When a (structured/unstructured) mesh is generated before hand, we could use a finite element/finite volume method. It is not easy to get high accuracy by using a finite volume method. The finite element method could have very high accuracy if high order basis functions are used. For elliptic boundary/interface problems, we could use Galerkin finite elements, the discontinuous Galerkin method, and the mixed finite element method. When the boundary/interface is complex, the apparent choice is to use a finite element method (FEM) with an unstructured mesh. However, it is not easy to generate an unstructured mesh especially when the boundary is very complex and the boundary changes with time. Another disadvantage of using FEM with an unstructured mesh is that it does not have the super convergence property which follows when using a uniform structured mesh. 

    Most of the comparison studies for elliptic boundary value/interface problems are conducted either through mesh refinement or by comparing methods using cartesian mesh \cite{LiuFedkiwKang2000, LeVequeLi1994, JohansenColella98, Mayo1984, MckenneyGreenGardMayo1995}. In this paper we are to perform a comparison study of two methods for solving the elliptic boundary value problem: the embedded boundary method using a cartesian mesh and the mixed finite element method using an unstructured mesh. The EBM uses a structured cartesian/rectangle grid. This method uses ghost cells along the boundary and the finite volume method to achieve 2rd accuracy in the potential and flux. The MFEM uses an unstructured triangular mesh. Instead of solving the second order elliptic equation, it solves two first order equations and gives the potential and flux at the same time. Higher order basis functions give higher order of accuracy. Refer to \cite{BrezziFortin1991} for a thorough discussion of mixed and hybrid finite element methods. For a more implementation oriented view, see \cite{Chavent1991}. The advantage and disadvantage of the MFEM are briefly discussed in \cite{Arnold1990}. For the comparison between FEM and MFEM, see the references cited in \cite{Chavent1999}. In this paper, we are to use the RT0 (Raviart-Thomas space of degree zero), the RT1 (Raviart-Thomas space of degree one), BDM1 (Brezzi-Douglas-Marini space of degree one) and BDM2 (Brezzi-Douglas-Marini space of degree two) as basis functions of the flux. We use the mixed-hybrid FEM. The final algebraic equations have only the potentials on the mesh edges as unknowns. To use the MFEM, we need to generate the mesh for the computational domain. There are mainly three methods for meshing: the Delaunay triangulation \cite{Shewchuk1997}, the advancing front method \cite{Marcum2001} and the quadtree/octree method \cite{YerryShephard1983}. In this paper, we use a method based on the quadtree/octree method. This method simplified the original construction by using marching cubes method to recover the interface. 

    The rest part of the paper is organized as follows. In section 2, we give the discretization of the two methods: embedded boundary method and the mixed finite element method. And also we will show briefly our method of generating the unstructured mesh for the mixed finite element. In section 3, we conduct the comparison study by solving a elliptic boundary problem with a known analytic solution. And in the last section, we give our conclusions. 

\section{Discretization}
     We are to solve the elliptic problem:
\begin{equation} 
\label{CS_ELLIPTIC}
\left\{ \begin{array}{l} \phi_{xx}+\phi_{yy} = f                 \\ \frac{\partial{\phi}}{\partial{n}} = g
\end{array} \right.
\end{equation}
in a complex domain, where $\phi(x,y)$ is called the potential. Since the gradient of the solution $\nabla\phi$ is often needed and more difficult to solve for, we will use the gradient errors as the comparison criterion. The gradients at both the regular grid centers and the boundary points are calculated and compared.

\subsection{Embedded Boundary Method}
The embedded boundary method is based on the finite volume discretization in grid blocks defined by the
rectangular Cartesian grid and the boundary. The solution is treated as a regular block centered quantity, even
when these centers are outside of the domain. However the gradient of the potential and the right hand side are
located in geometrical centers (centroids) of the partial grid blocks cut by the boundary \cite{JohansenColella98}.
This treatment has advantages when dealing with geometrically complex domains; it also ensures
second-order accuracy of the solution. \\

In the 2D case, each regular grid block is a square. Using the divergence theorem and
integrating the flux  $ {\mathbf F}=\nabla\varphi $ over the control volume, the differential operator can be
discretized as
\begin{equation}
\label{discretization} (L\varphi)_{\Delta_i}=\frac{1}{V_i}(\sum _j {{\mathbf F}_j \cdot {\mathbf S_j}}),
\end{equation}
where $V$ and ${\mathbf S}$ are size of the control volume and block edge respectively, and ${\mathbf F}$ is
the flux across the geometric center of each edge. For full edges (not cut by the boundary),
${\mathbf F}_j$ is obtained by the central difference while the flux across the partial block edges is obtained
using a linear interpolation between centered difference fluxes in adjacent blocks.

\begin{figure}
\centering
\includegraphics[scale=0.5]{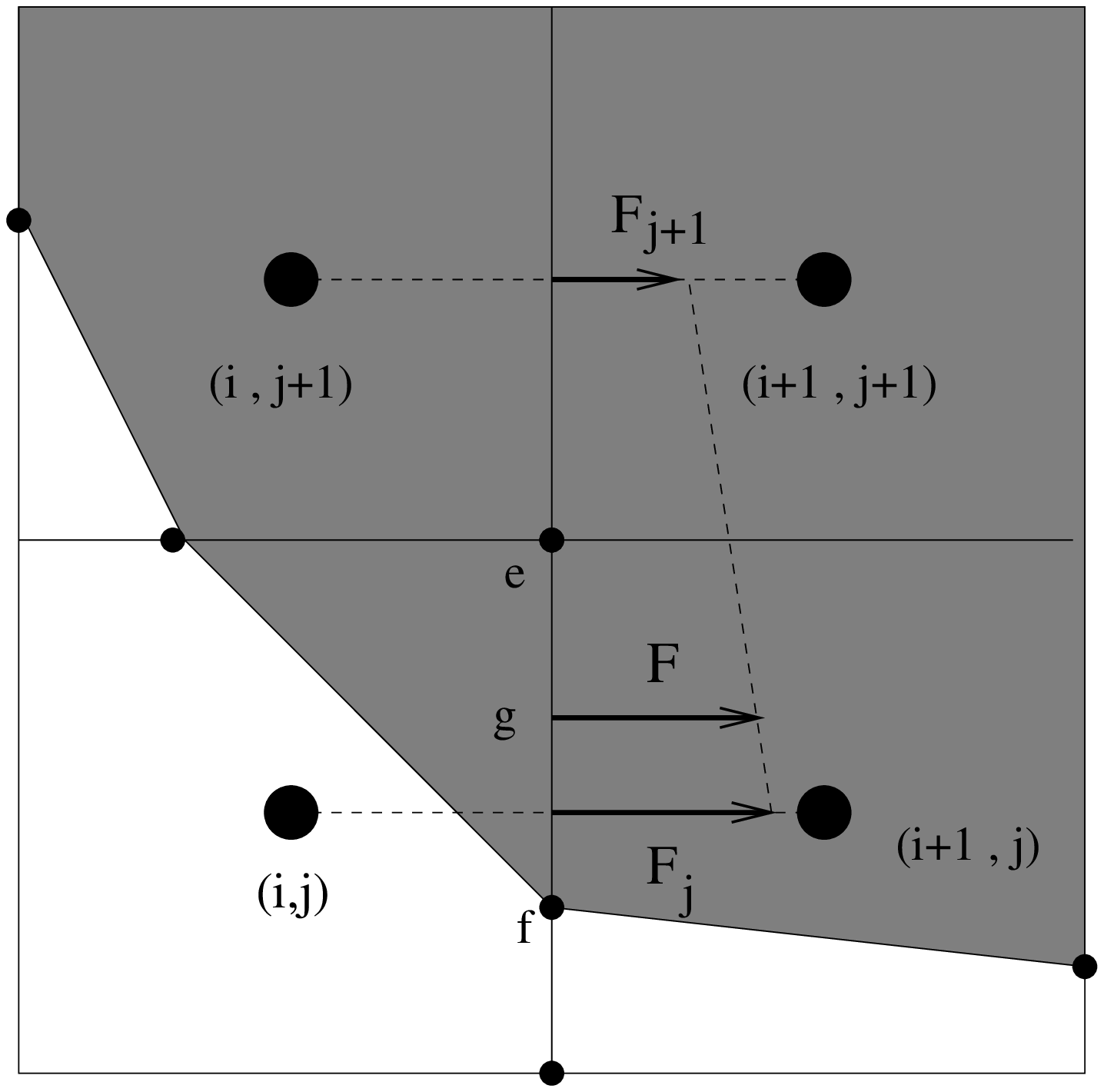}%
\hspace{1in}%
\includegraphics[scale=0.5]{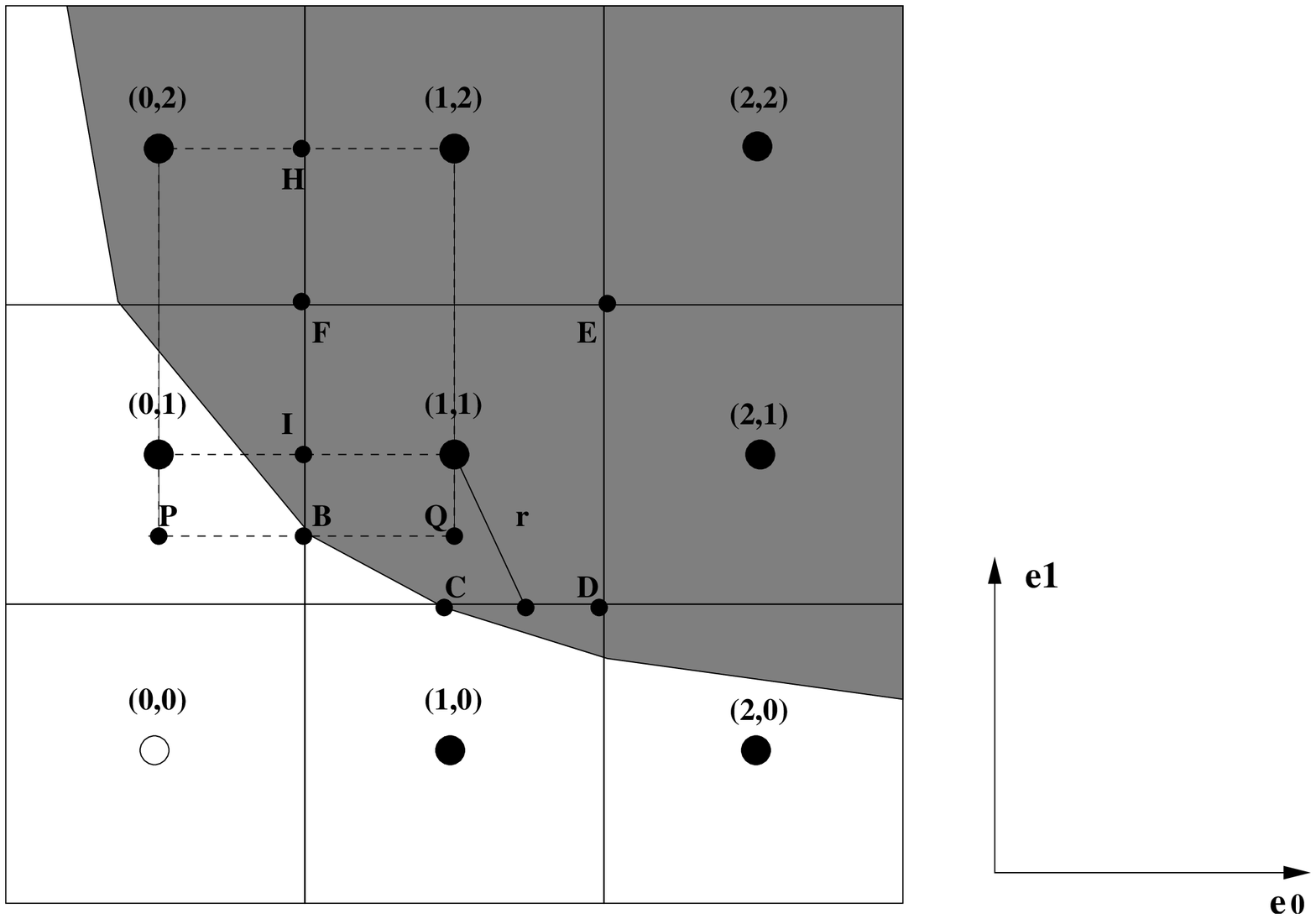}
\caption{Linear interpolation of flux(Left) and Stencil setting(Right)} \label{stencil}
\end{figure}

The flux interpolation method is illustrated as the left side of Fig. \ref{stencil}. The flux across the center
$g$ of the partial edge $ef$ is obtained using the linear interpolation between the fluxes ${\mathbf F}_j$ and
${\mathbf F}_{j+1}$, which are the finite differences of potentials at the centers of the corresponding regular grid
blocks. The flux at the domain boundary is given by the Neumann condition.

In order to implement the embedded boundary method, the boundary is reconstructed using its intersections with
grid lines. The following assumptions and simplifications are made:
\begin{flushleft}
\begin{enumerate}
\item The maximum number of intersection of each block edge with the boundary curve is one.
\item The elliptic problem domain within each grid block forms a connected set.
\item The positions of the boundary points are adjusted to remove partial blocks with volumes less than a certain preset value. 
\end{enumerate}
\end{flushleft}

The first and second assumptions are generally satisfied when the curvature of the boundary curve is not too large or
the mesh is sufficiently refined. The third one is necessary since blocks of arbitrary small volumes introduce
large numerical errors and increase the condition number of the linear system resulting from the discretization.

The summary of the algorithm is as follows.

(1) The elliptic domain boundary is constructed using intersection points of the grid free boundary with grid
lines. All the grid blocks are divided into three types: INTERNAL, PARTIAL, and
EXTERNAL, which means completely within, partially within (cut by the boundary), and completely outside of
the elliptic domain.

(2) The number of blocks marked as PARTIAL or INTERNAL is counted and the total size of the linear system is
set. For each block marked as PARTIAL, all block edges are also divided into three types similar to the types introduced above. The
center position and length of each partial edge are stored. A 9-point stencil is set to calculate fluxes across the
control volume $BADEF$, as shown in the right side of Figure \ref{stencil}, where the elliptic problem domain is
the shaded region and the filled circles represent locations where the potential is defined. We define a $3\times3$ matrix
$C$ with matrix elements $c(i,j)$ representing the coefficient of $\varphi$ centered at $(i,j)\, (i,j=0,1,2)$
according to ${\mathbf F}=\sum_{i,j}c(i,j)\varphi(i,j)$. Therefore $\varphi (1,1)$ is always the potential
located within the control volume. We further denote $c(i,j)$ as $c({\mathbf V})$, where the vector ${\mathbf V}$
has components $i$ and $j$. The vector ${\mathbf r}$ is drawn from the regular block center to the center of the
block edge on which the flux is to be integrated. Suppose the basis of the Cartesian coordinate is formed by
the unit vectors ${\mathbf e}_i (i=0,1)$, and ${\mathbf e}$ is the vector with all unit elements, then the unit
vector ${\mathbf e'}_i=\textrm{sign}({\mathbf r} \cdot {\mathbf e}_i){\mathbf e}_i$ gives orientational information of
${\mathbf r}$. For the linearly interpolated flux $F$ along a direction $d$, the corresponding coefficients
are:

$$c({\mathbf e})=\frac{a-1}{h_d}; \quad c({\mathbf e} +{\mathbf e'}_d)=\frac{1-a}{h_d}$$

$$c({\mathbf e} + {\mathbf e'}_{d'})=\frac{-a}{h_d}; \quad c({\mathbf e} + {\mathbf e'}_d + {\mathbf e'}_{d'})=\frac{a}{h_d}$$
where $d',d$ = 0,1 and $d' \neq d$, $h_d$ is the grid spacing in the direction $d$, and $a=\frac{|{{\mathbf r}}
\cdot {{\mathbf e}_{d'}}|}{h_{d'}}$ is the block edge aperture.

(3) Substituting ${\mathbf F}=\sum_{i,j}c(i,j)\varphi(i,j)$ into the equation (\ref{discretization}) and summing fluxes through all edges of each PARTIAL block, the coefficients at stencil points are set and added to the
global matrix. Note that the right hand side in equation (\ref{discretization}) must be evaluated at the centroid of
the partial block.
\\
(4) The resulting linear system $A{\mathbf x}={\mathbf b}$ is solved. Then the gradient of the potential is
calculated at all PARTIAL and INTERNAL block centers, even if these centers are outside of the elliptic domain.
Either the centered difference or quadric interpolation is used to maintain the second order accuracy. For
example, the x-derivative of the potential $\varphi_x(1,1)$ can be
easily calculated by the centered differences of  $\varphi(0,1)$ and $\varphi(2,1)$. For some points
near the boundary, such as point located at (0,2), $\varphi_x(0,2)$ is obtained as the $x$ derivative of the
quadric curve, which interpolates potential values $\varphi(0,2)$, $\varphi(1,2)$, and  $\varphi(2,2)$. We also
calculate the gradient of the potential at the boundary points for later comparison. Take the point $B$ in  Fig. \ref{stencil} for example, ${\varphi}_x(B)$ is the extrapolation between ${\varphi}_x(H)$ and ${\varphi}_x(I)$, which
in turn are calculated by centered differences. To calculate ${\varphi}_y(B)$, first we extrapolate
${\varphi}_y(P)$ from ${\varphi}_y(0,1)$ and ${\varphi}_y(0,2)$,${\varphi}_y(Q)$ from
${\varphi}_y(1,1)$ and ${\varphi}_y(1,2)$. Then ${\varphi}_y(B)=\frac{1}{2}({\varphi}_y(P)+{\varphi}_y(Q))$.

\subsection{Mixed Finite Element Method}

     The mixed finite element method(MFEM) solves for the potential $\phi$ and the flux $\nabla\phi$ at the same time. Thus, it solves
\begin{displaymath}    
    \vec{q} = -a\nabla \phi
\end{displaymath}
\begin{displaymath}
    \nabla \cdot \vec{q} = f
\end{displaymath}
For the mixed finite element method, two function spaces are needed: one scalar space for the potential $\phi$ and one vector space for the flux $\vec{q}$. The unknowns are potentials on the elements and flux on the edges. To reduce the problem to a smaller one, the mixed-hybrid finite element is modified by introducing a Lagrangian multiplier on the edges. Chavent and Roberts \cite{Chavent1991} give in detail an implementation using rectangle elements. The final algebraic equations only have TPs (the Lagrangian multiplier, also the potential on the edges) as unknowns, thus reducing the number of unknowns. Later they reduced the problem further by introducing an unknown variable defined inside the element \cite{Chavent1999}. Now instead of unknowns defined on edges, they have only one unknown in each triangle element. Since the number of triangles is much smaller than the number of edges, the problem is reduced into a smaller one. However they only use the lowest order RT basis in their derivation. Since we do not know whether their approach could be extended to use higher order basis functions, our implementation uses the first approach \cite{Chavent1991}. See \cite{BrezziFortin1991} for a more theoretical treatment of the subject.

    In this paper, we use the mixed-hybrid finite element with four different basis functions for the flux: the RT0 (Raviart-Thomas space of degree zero), the RT1 (Raviart-Thomas space of degree one) and BDM1 (Brezzi-Douglas-Marini space of degree one) and BDM2. 
     The basis function for the flux in RT0 is:
\begin{displaymath}
    \vec{s}|_{K} =    a \left(\begin{array}{l}x\\y\end{array} \right) + \left(\begin{array}{l}b\\c\end{array} \right),
\end{displaymath}
The basis function for the potential is a constant. The lagrangian multiplier TP defined on edges is also constant. RT0 has 1st order accuracy in the flux and potential in the $L_2$ norm.

     The basis function for the flux in BDM1 is
\begin{displaymath}
    \vec{s}|_{K} =    \left(\begin{array}{l} a_1x+ a_2y+ a_3 \\ b_1x+b_2y+b_3\end{array} \right),
\end{displaymath}
The potential is also constant. However TP is linear. The accuracy for BDM1 is 2nd order in the flux and potential.

     The basis function for the flux in RT1 is
\begin{displaymath}
    \vec{s}|_{K} =    \left(\begin{array}{l} a_1x+ a_2y+ a_3 \\ b_1x+b_2y+b_3\end{array} + \begin{array}{l} x\\ y\end{array} \times (c_1x+c_2y)\right),
\end{displaymath}
and the basis for the potential is a linear function:
\begin{displaymath}
    v|_{K} = P_1 p_1 +P_2 p_2 +P_3 p_3, 
\end{displaymath} 
where $p_1, p_2, p_3$ are the basis functions. TP is also linear. The accuracy is the same as BDM1.

    The basis function for the flux in BDM2 is
\begin{displaymath}
    \vec{s}|_{K} =    \left(
                        \begin{array}{l} 
                            a_1x^2 + a_2xy + a_3y^2 + a_4x + a_5y + a_6 \\
                            b_1x^2 + b_2xy + b_3y^2 + b_4x + b_5y + b_6
                        \end{array} 
                      \right),
\end{displaymath}
The potential is the same as for RT1. TP is quadratic. The accuracy for BDM2 is 3rd order in the flux and potential.    

For implementation in details, see \cite{Wang2006, Chavent1991}.

\subsection{Mesh Generation}
Here, we first introduce our mesh generation method briefly and then give the point location algorithm to locate the triangle which contains a given point. Refer to \cite{Wang2006} for more detail.
\subsubsection{Quadtree Mesh Generation}
     For the mixed finite element method, we use an unstructured mesh with triangles only. Our method for mesh generation is similar to the Quadtree/Octree based mesh generation method developed by Yerry and Shephard \cite{YerryShephard1983}. However, there is one important simplification in the interface recovering step. The quadtree/octree is a tree structure \cite{Samet1990}. Each quadrant in the quadtree has exactly four children and each octant in the octree has exactly eight children. The quadtree/octree is used widely and it is used here for automatic mesh refinement (AMR). The quadtree/octree data structure has a number called level, representing the depth of the tree structure. The root has level 0, its four children has level 1 and so on.

The quadtree/octree mesh generation method is simple and it consists of the following steps (using quadtree as example):
\begin{flushleft}
\begin{enumerate}
\item Partition the computational region into a quadtree with the level difference between neighbor quadrants being at most 1. Now all those quadrants are either full interior quadrants or partial/boundary quadrants.
\item Triangulate the full interior quadrants.
\item Triangulate the partial quadrants to recover the interface.
\item Post processing the mesh. If we used templates to triangulate the partial quadrants and recover the interface in step 3, we need to move those interface points onto the interface in the post processing step.
\end{enumerate}
\end{flushleft}

The main difference of our method compared with Yerry and Shephard's method \cite{YerryShephard1983} lies in the 3rd step in recovering the interface. In their original method, the interface could cross over the edges and on the vertices of the quadtree. In our method, we assume that the interface could only cross over the edges and for each edge, there is at most one crossing. Thus we are using the the marching cubes method (MC) for interface recovering. The marching cubes method was proposed by Lorensen and Cline \cite{LorensenCline1987} for extracting an isosurface from volumetric data. Here we use it to recover our interface. Note that if the input uses a boundary representation (edges), then our mesh might not be conforming to the input boundary. The reason is apparent when we check the way the constraint delaunay triangulation method recovers the interface: We need to recover the vertices first and then the edges in 2D. For 3D problem, we need to recover the vertices, edges and faces. In our method we recover only edges in 2D and only faces in 3D.

\subsubsection{Quadtree Mesh Point Location}
After the mesh is given, we use a finite element to set up the matrix and solve for the unknowns. Sometimes we need the solutions for an arbitrary point inside the mesh, which is in fact a point location problem: find the triangle/tetrahedron which contains the given point. The point location problem and another closely related problem called the range search problem are two famous problems in computational geometry. See \cite{Berg1998} and references cited therein.

If only one point is queried, we only need to loop through every triangle/tetrahedron of our mesh and test whether the triangle/tetra contains the given point. The time complexity is clearly $O(N)$ where $N$ is the number of triangles/tetrahedra inside the mesh. If $m$ such points are to be queried, such an approach would not be applicable when $m$ is large such as $m=O(N)$. We would be in such a situation if we solve an elliptic interface problem using the mixed finite element on an unstructured grid and then interpolate the flux back onto an cartesian grid.

To speed up the point location problem, it is a common practice to preprocess the mesh and set up some special data structure. Fortunately, we do not need to create a new data structure here. Since the quadtree/octree is a tree structure, we use it for the point location.
Our algorithm is the following: \\
Given point P, the Quadtree/Octree and mesh,
\begin{flushleft}
\begin{enumerate}
\item first use the quadtree/octree structure to find a leaf quadrant/octant;
\item second use the leaf quadrant/octant to find an triangle/tetrahedron which would be used as an starting point to find the target triangle/tetrahedron;
\item walk through the mesh to the given point P.
\end{enumerate}
\end{flushleft} 
\section{A Comparison Study}

\begin{figure}[htp]
\centering
\includegraphics[scale=0.3]{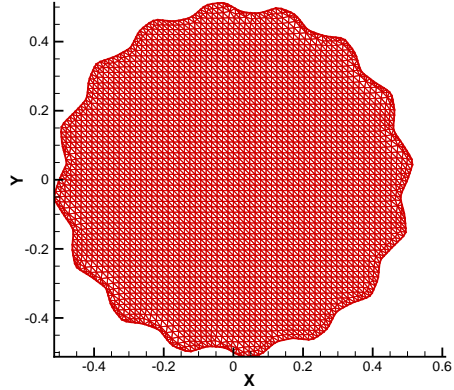}
\caption{The unstructured computational mesh for a $128\times128$ mesh}
\label{Mesh_7x7}
\end{figure}

For our test problem, we use $\phi = e^{\frac{x^2+y^2}{2}}$ as the exact solution of the elliptic equation (\ref{CS_ELLIPTIC}); f and g are obtained by differentiating $\phi$. We will show two different testing problems using the same equation and analytic solution. The difference between the two problems lies only in the different boundaries: the second boundary is more complex than the first one. The EBM uses a structured cartesian grid. The mixed finite element methods use an unstructured grid based on the quadtree/octree construction. The quadtree/octree have minimum and maximum levels. In order to compare the results, we need to have comparable grids by letting the minimum/maximum level of the quadtree to be equal. Fig. \ref{Mesh_7x7} shows the grid used by the mixed finite element methods when EBM uses the $128\times128$ grid. Thus the mesh is uniform. We compare the results using the $L_2$ norm of the flux $\|\nabla\phi\|_2$. The norm is defined as:
\begin{equation}
\label{CS_DISCRETE_NORM}
\left\{
\begin{array}{l}
    \|\nabla\phi\|_2 = \sqrt{\sum_{\textrm{face}}\|\nabla\phi\|_{2,\textrm{face}}^2}      \\
    \|\nabla\phi\|_{2,\textrm{face}} = \textrm{Area}(\textrm{face}) \times \sqrt{\phi_x(x_0,y_0)^2+\phi_y(x_0,y_0)^2}
\end{array}
\right.
\end{equation}
where $(x_0,y_0)$ is the center of the rectangle for the cartesian grid used by the EBM. For the MFEM, we first interpolate the fluxes at the center of the cartesian grid, and then compute the norm.

    The matrices for both methods are solved using methods in the PETSc package. Here we use the BiCGSTAB method with the ilu method as preconditioner. We have tried different methods (such as lu, Cholesky, CG, GMRES, BiCGSTAB etc) with different preconditioners in the PETSc packages and find that the BiCGSTAD method with ilu as preconditioner is the fastest for solving our matrices.

\subsection{Embedded Boundary Method vs. Mixed Finite Element Method}

\begin{figure}[htp]
\centering
\includegraphics[scale=0.3]{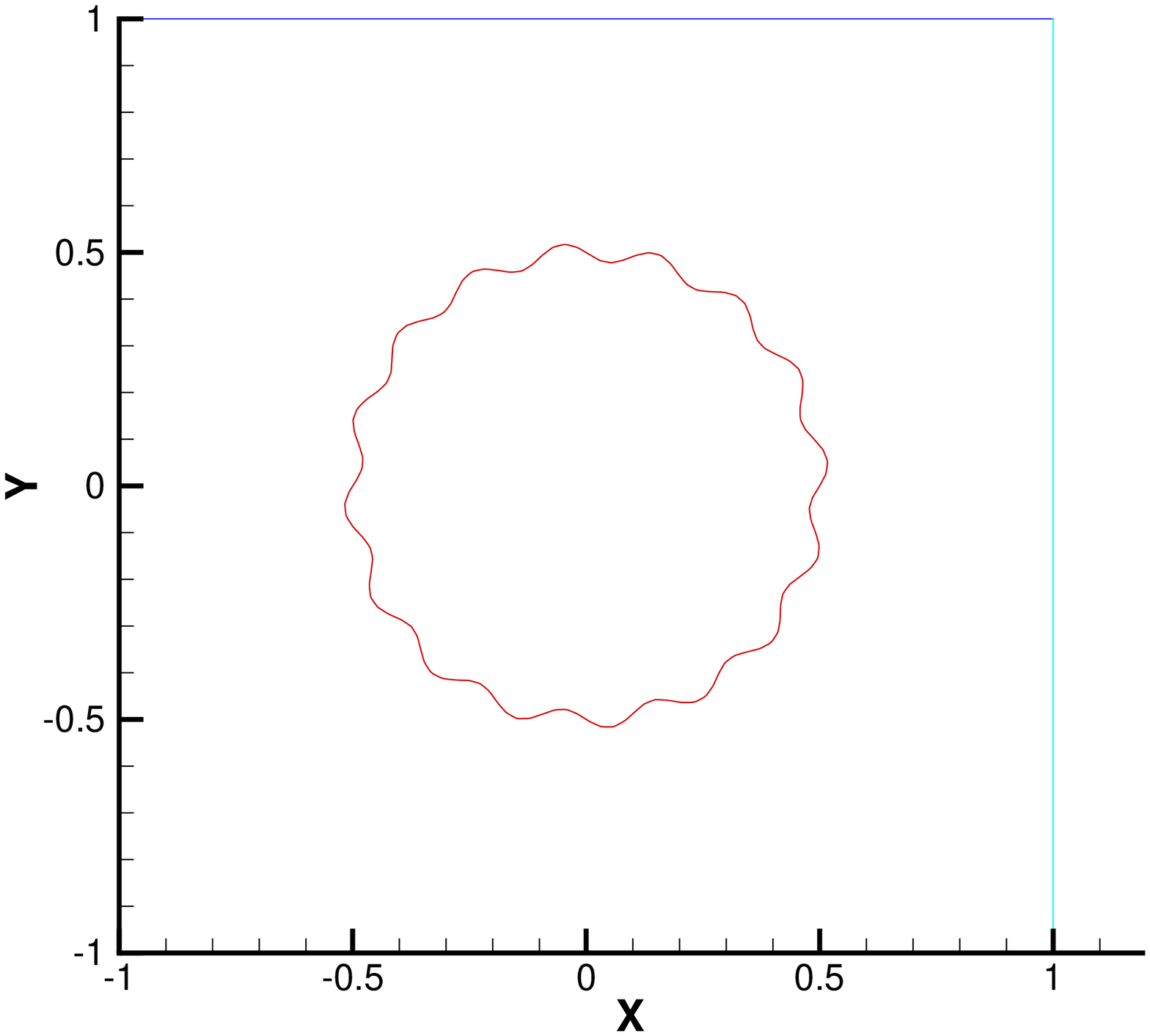}
\caption{Boundary for the first test}
\label{INTERFACE0}
\end{figure}

The first problem uses a simple boundary. The computational domain lies inside a perturbed circle as in Fig. \ref{INTERFACE0}.

\begin{table}
\caption{Convergence and Timing Study using Uniform Mesh for the Boundary in Fig. \ref{INTERFACE0}}
\label{INTERFACE0_TABLE_CONV_TIMING}
\begin{tabular}{|c||r|r|r|r|r|}
\hline
Mesh            &\multicolumn{5}{c|}{EBM} \\
\cline{2-6}
Size            &  error        &  ratio    &  time     &iterations &unknowns\\
\hline
$64\times64$    &1.593569e-04   & N/A       &0.016966   & 32        &861   \\
$128\times128$  &3.670301e-05   & 2.118     &0.099232   & 60        &3338   \\
$256\times256$  &8.686625e-06   & 2.099     &0.699156   & 116       &13160   \\
$512\times512$  &2.134996e-06   & 2.074     &6.023590   & 242       &52056   \\
\hline
\end{tabular}

\begin{tabular}{|c||r|r|r|r|r|}
\hline
Mesh            &\multicolumn{5}{c|}{RT0} \\
\cline{2-6}
Size            &  error        &  ratio    &  time     &iterations &unknowns\\
\hline
$64\times64$    &1.011978e-03   & N/A       &0.193374   & 74        & 2642\\
$128\times128$  &5.121661e-04   & 0.982     &0.836596   & 108       & 10141\\
$256\times256$  &2.651845e-04   & 0.950     &5.723007   & 219       & 39751\\
$512\times512$  &1.353009e-04   & 0.971     &42.336410  & 462       & 156715\\
\hline
\end{tabular}

\begin{tabular}{|c||r|r|r|r|r|}
\hline
Mesh            &\multicolumn{5}{c|}{BDM1}    \\
\cline{2-6}
Size            &  error        &  ratio    &  time     &iterations &unknowns\\
\hline
$64\times64$    &1.329723e-04   & N/A       &0.475765   & 87        & 5286\\
$128\times128$  &3.715907e-05   & 1.839     &3.131988   & 171       & 20284\\
$256\times256$  &9.794951e-06   & 1.924     &19.277468  & 306       & 79504\\
$512\times512$  &2.538575e-06   & 1.948     &141.012044 & 597       & 313432\\
\hline
\end{tabular}

\begin{tabular}{|c||r|r|r|r|r|}
\hline
Mesh            &\multicolumn{5}{c|}{RT1} \\
\cline{2-6}
Size            &  error    &  ratio    &  time     &iterations &unknowns\\
\hline
$64\times64$    &1.701312e-05   & N/A       &0.807441   & 87        &   5286\\
$128\times128$  &4.628462e-06   & 1.878     &4.472460   & 172       &  20284\\
$256\times256$  &1.215990e-06   & 1.928     &24.581123  & 305       &  79504\\
$512\times512$  &3.125208e-07   & 1.960     &163.794142 & 607       & 313432\\
\hline
\end{tabular}

\begin{tabular}{|c||r|r|r|r|r|}
\hline
Mesh            &\multicolumn{5}{c|}{BDM2} \\
\cline{2-6}
Size            &  error    &  ratio    &  time     &iterations &unknowns\\
\hline
$64\times64$    &4.598191e-07   & N/A       &1.242945   & 100       &  7929\\
$128\times128$  &4.169450e-08   & 3.463     &6.773787   & 191       & 30426 \\
$256\times256$  &4.824547e-09   & 3.111     &40.007515  & 317       &119256 \\
$512\times512$  &2.865473e-09   & 0.751     &336.560434 & 756       &470148 \\
\hline
\end{tabular}
\end{table}

Table \ref{INTERFACE0_TABLE_CONV_TIMING} displays the errors and timing results for different mesh sizes. The convergence ratios with mesh refinement, the number of unknowns for the linear system, and the number of iterations for the linear solver are also listed. The maximum relative tolerance is $1e^{-9}$. The errors are measured by the $L_2$ norm of $\nabla\varphi$ defined by (\ref{CS_DISCRETE_NORM}). From the table, we see that RT0 has only first order accuracy, EBM/BDM1/RT1 have 2nd order accuracy and BDM2 has 3rd order accuracy. The EBM is much faster than the other four methods when the same mesh size is used. The most apparent reason is that it has fewer unknowns than the other four methods. As expected, the RT0 is faster than BDM1/RT1/BDM2 since it has at most one half the number of the unknown variables. However, RT0 only has 1st order accuracy. Although BDM1/RT1 are both 2nd accurate method with the same number of unknowns, they have different characteristics in their timing and accuracy. The BDM1 is less accurate but faster for given mesh size. BDM2 has the highest order accuracy of all five methods. For the same order of accuracy, the fastest method is BDM2, then EBM/RT1/BDM1/RT0.

Table \ref{INTERFACE0_RT0_TIMING_DETAIL} gives the timing of the RT0 method for the mesh generation, the matrix setup/solve and the interpolation of the solution. Note that RT0/BDM1/RT1/BDM2 use the same mesh. Therefore their mesh generation time is the same. Their timing differences lie only in the matrix setup/solve step. Here we find that the time spent on generating the mesh is only a small part of the total time when the mesh size is large. Most of the time are spent on solving the algebraic
equation (timing for the matrix setup is comparable with that of the interpolation step). It is more apparent when the mesh size is increased. For example, the ratio of time spent on the matrix setup/solve step compared to the mesh generation step is about 1.21 when the $64\times64$ mesh is used. The same ratio increases to 9.85 when the $512\times512$ mesh is used.
\begin{table}
\caption{timing of RT0 in detail}
\label{INTERFACE0_RT0_TIMING_DETAIL}
\begin{tabular}{|c||l|l|l|}
\hline
Mesh            &\multicolumn{3}{c|}{RT0} \\
\cline{2-4}
 Size           &  mesh     & matrix setup/solve    & interpolation\\
\hline
$64\times64$    & 0.083814  & 0.101778              & 0.004258  \\
$128\times128$  & 0.257736  & 0.547282              & 0.016645  \\
$256\times256$  & 0.983702  & 4.613645              & 0.073030  \\
$512\times512$  & 3.849851  & 37.933463             & 0.343644  \\
\hline
\end{tabular}
\end{table}

\begin{figure}[htp]
\centering
\includegraphics[scale=0.3]{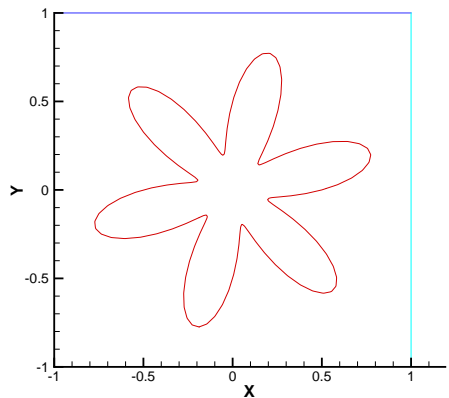}%
\caption{The boundary for the second test}
\label{INTERFACE1}
\end{figure}

In the following, we use the EBM and MFEM to solve the same problem but using a more complicated boundary as shown in Fig. \ref{INTERFACE1}. The errors are still measured by the $L_2$ norm of $\nabla\varphi$ defined by (\ref{CS_DISCRETE_NORM}) and the max tolerance is $1e^{-9}$. The mesh is more refined in order to well resolve the boundary. Table \ref{INTERFACE1_TABLE_CONV_TIMING} shows the convergence and timing results of the five methods. The general conclusion is the same as the first test. The RT0 is 1st order accurate, the EBM/BDM1/RT1 method are 2nd order and BDM2 is 3rd order accurate in flux. The EBM method is still the fastest method for the same mesh size. For the same accuracy, we have BDM2, then EBM/RT1/BDM1/RT0 in decreasing order of speed. 

\begin{table}
\caption{Convergence and Timing Study using Uniform Mesh for the Boundary in Fig. \ref{INTERFACE1}}
\label{INTERFACE1_TABLE_CONV_TIMING}
\begin{tabular}{|c||r|r|r|r|r|}
\hline
Mesh          &\multicolumn{5}{c|}{EBM}   \\
\cline{2-6}
 Size           &  error        &  ratio &  time    &iterations & unknowns\\
\hline
$64\times64$    &2.110753e-04   & N/A    &0.022319  & 43        &1008    \\
$128\times128$  &5.779287e-05   & 1.869  &0.164115  & 91        &4008   \\
$256\times256$  &1.472989e-05   & 1.920  &1.438516  & 209       &15738   \\
$512\times512$  &3.641386e-06   & 1.952  &10.398294 & 360       &61967  \\
\hline
\end{tabular}

\begin{tabular}{|c||r|r|r|r|r|}
\hline
Mesh            &\multicolumn{5}{c}{RT0} \\
\cline{2-6}
 Size           &  error        &  ratio    &  time     &iterations &unknowns\\
\hline
$64\times64$    & 1.806532e-03  & N/A       &0.283352   & 115       & 3177   \\
$128\times128$  & 1.142133e-03  & 0.661     &1.624428   & 218       & 12341  \\
$256\times256$  & 6.140341e-04  & 0.895     &11.236136  & 415       & 47870  \\
$512\times512$  & 3.166839e-04  & 0.955     &79.363582  & 770       & 187229 \\
\hline
\end{tabular}

\begin{tabular}{|c||r|r|r|r|r|}
\hline
Mesh            &\multicolumn{5}{c}{BDM1} \\
\cline{2-6}
 Size           &  error        &  ratio    &  time     &iterations &unknowns\\
\hline
$64\times64$    & 1.695040e-04  & N/A       &0.832989   & 151       & 6354   \\
$128\times128$  & 5.857866e-05  & 1.533     &5.360611   & 278       & 24682   \\
$256\times256$  & 1.641488e-05  & 1.835     &33.627868  & 461       & 95740  \\
$512\times512$  & 4.311759e-06  & 1.929     &320.587307 & 1185      & 374458  \\
\hline
\end{tabular}

\begin{tabular}{|c||r|r|r|r|r|}
\hline
Mesh            &\multicolumn{5}{c}{RT1} \\
\cline{2-6}
 Size           &  error        &  ratio    &  time     &iterations &unknowns\\
\hline
$64\times64$    &2.203143e-05   & N/A       &1.212332   & 151       & 6354\\
$128\times128$  &7.323467e-06   & 1.589     &7.185103   & 296       & 24682\\
$256\times256$  &2.032626e-06   & 1.849     &45.390103  & 549       & 95740\\
$512\times512$  &5.312876e-07   & 1.936     &312.086512 & 1055      & 374458\\
\hline
\end{tabular}

\begin{tabular}{|c||r|r|r|r|r|}
\hline
Mesh            &\multicolumn{5}{c}{BDM2} \\
\cline{2-6}
 Size           &  error        &  ratio    &  time     &iterations &unknowns\\
\hline
$64\times64$    &6.651279e-07   & N/A       &1.916726   & 176       & 9531\\
$128\times128$  &7.259815e-08   & 3.196     &12.185872   & 323       & 37023\\
$256\times256$  &9.293894e-09   & 2.966     &90.869104  & 660       & 143610\\
$512\times512$  &3.087661e-09   & 1.590     &607.406103 & 1162      & 561687\\
\hline
\end{tabular}
\end{table}

\begin{figure}[htp]
\centering
\includegraphics[scale=0.3]{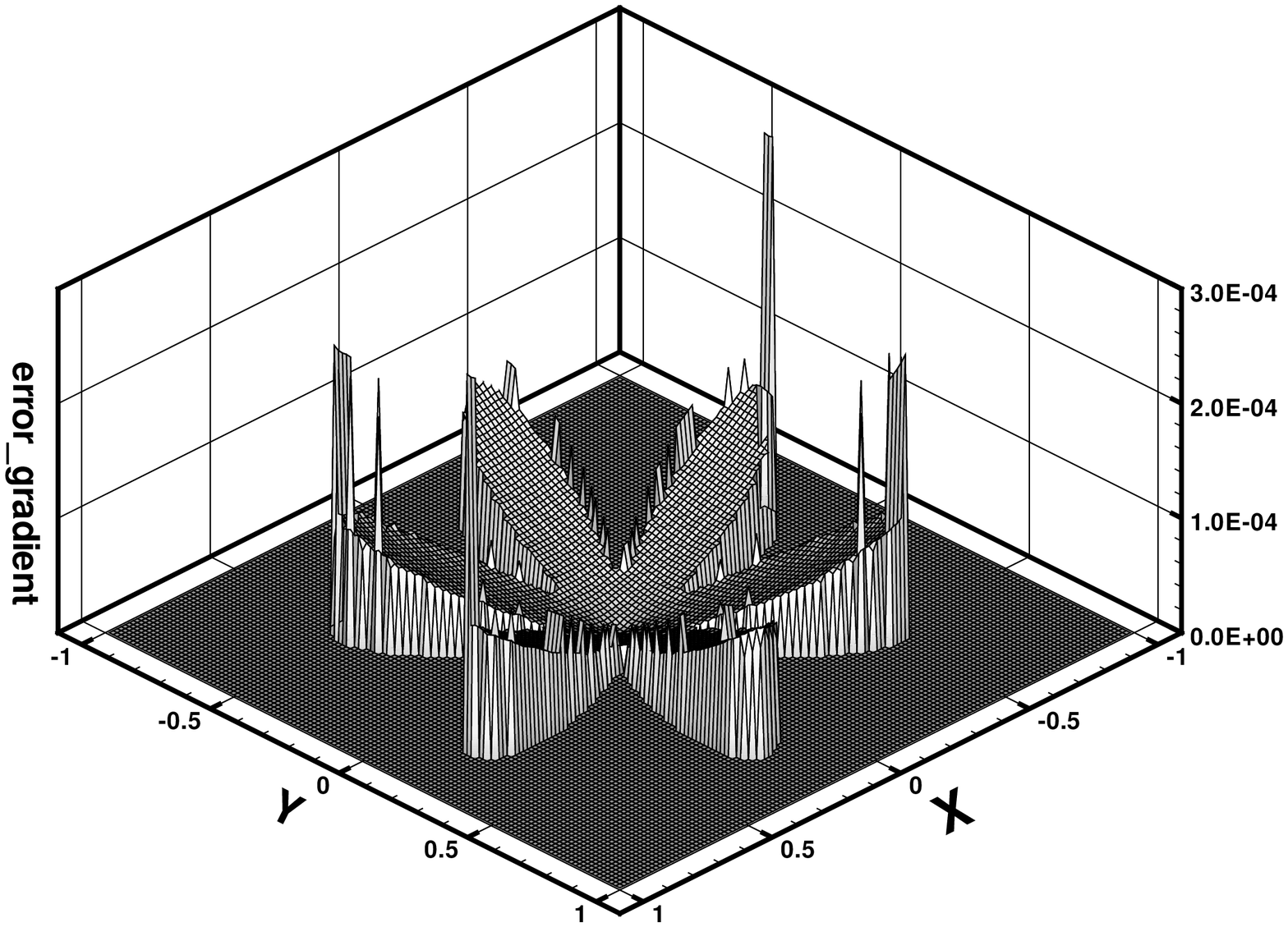}
\caption{norm of gradient error by EBM using the $128\times128$ grid}\label{EBM_DPhi_Error}
\end{figure}

\begin{figure}[htp]
\centering
\includegraphics[scale=0.3]{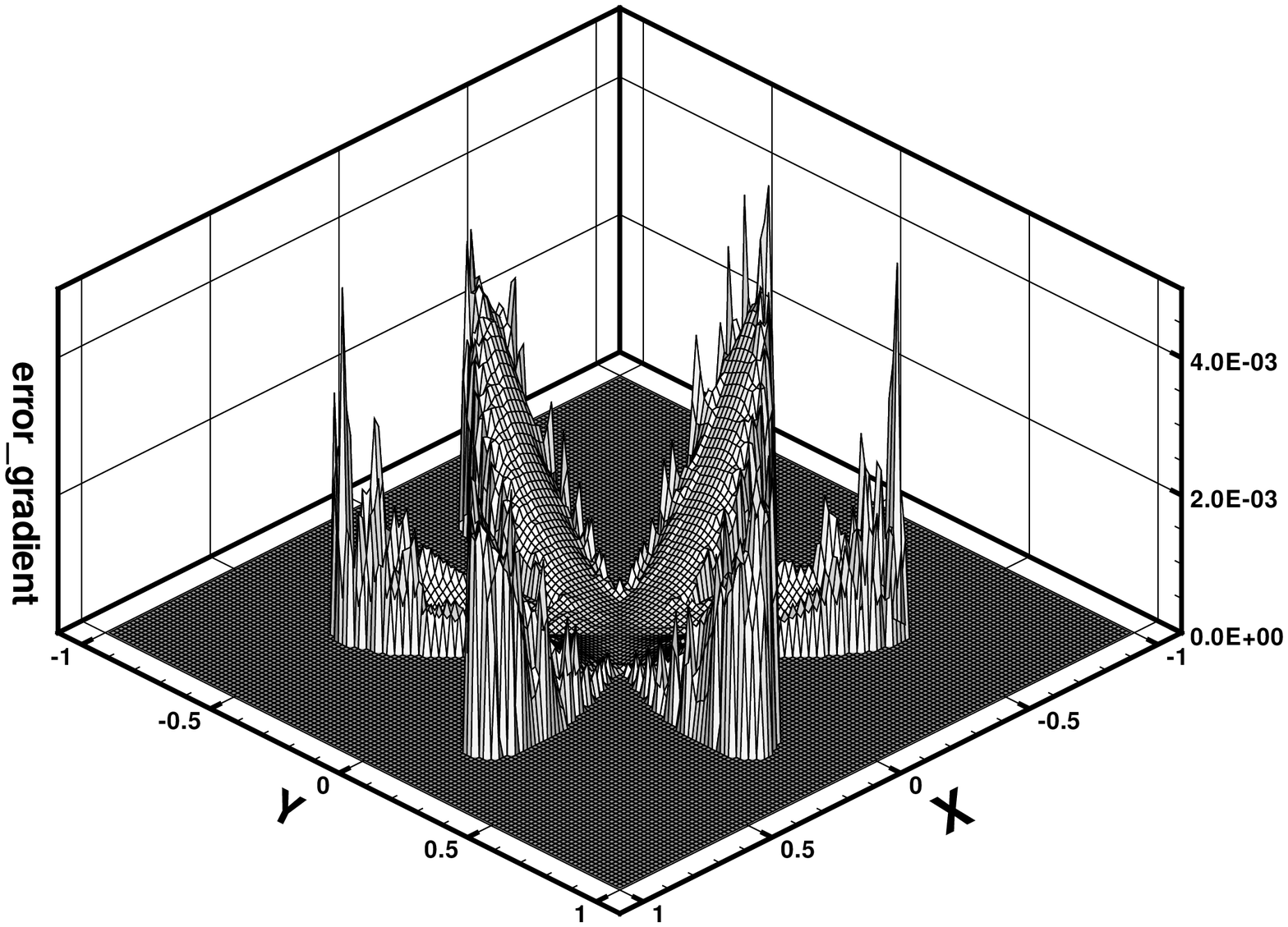}
\caption{norm of gradient error by RT0 using the $128\times128$ grid}\label{RT0_DPhi_Error}
\end{figure}

\begin{figure}[htp]
\centering
\includegraphics[scale=0.3]{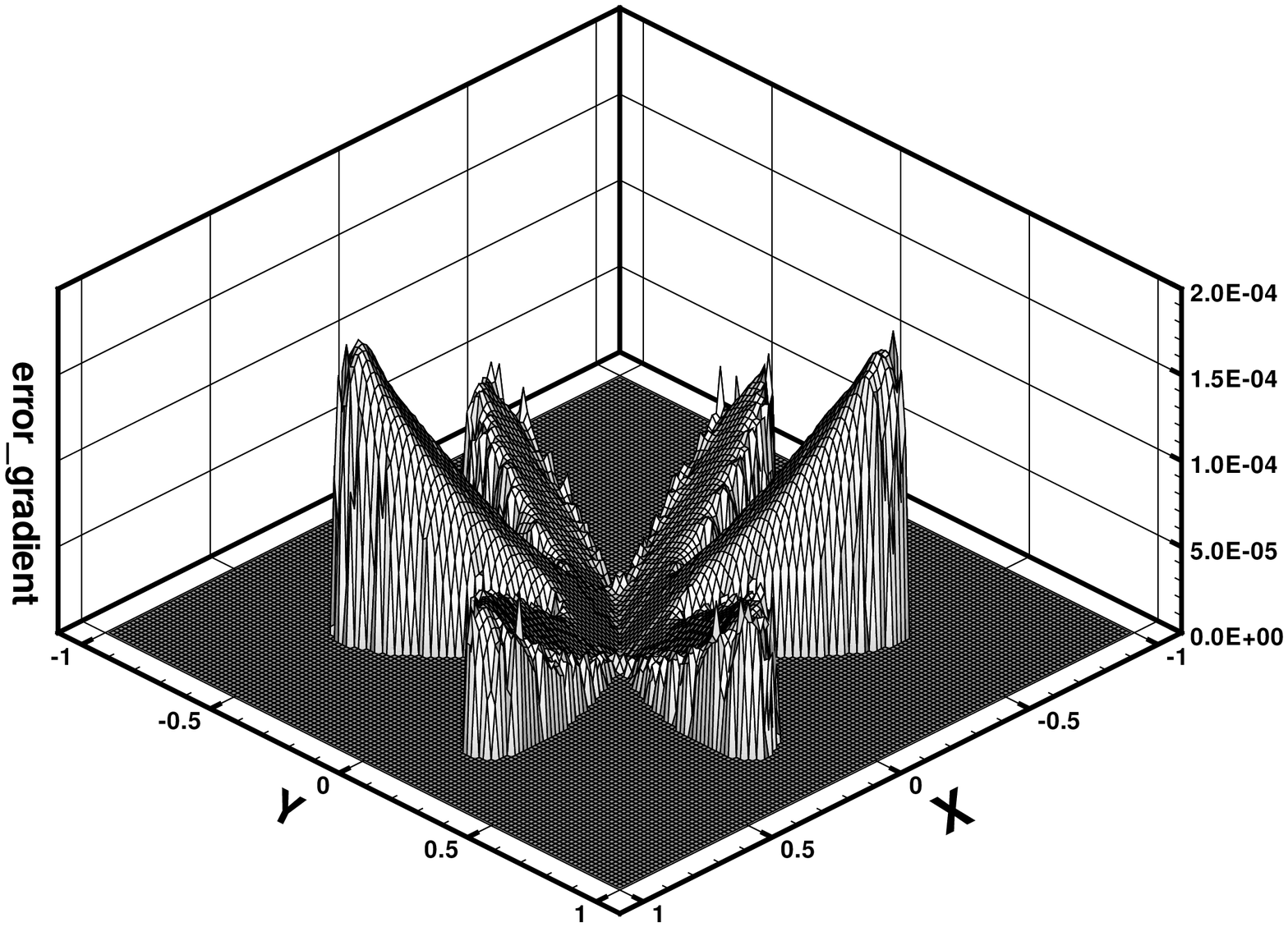}
\caption{norm of gradient error by BDM1 using the $128\times128$ grid}\label{BDM1_DPhi_Error}
\end{figure}

\begin{figure}[htp]
\centering
\includegraphics[scale=0.3]{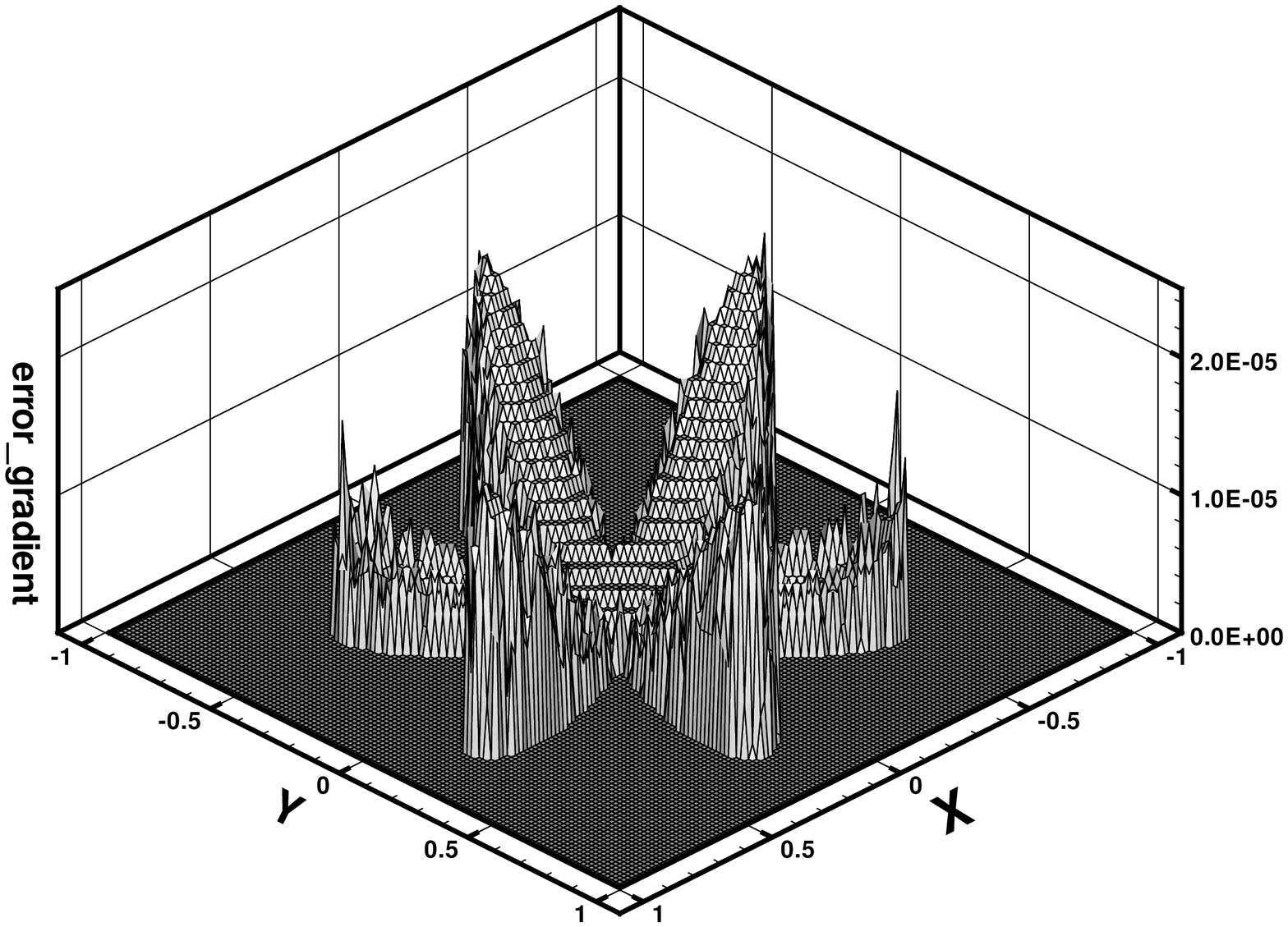}
\caption{norm of gradient error by RT1 using the $128\times128$ grid}\label{RT1_DPhi_Error}
\end{figure}

\begin{figure}[htp]
\centering
\includegraphics[scale=0.3]{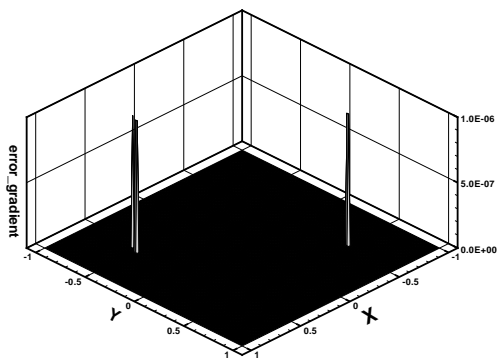}
\caption{norm of gradient error by BDM2 using the $128\times128$ grid}\label{BDM2_DPhi_Error}
\end{figure}

Figs. \ref{EBM_DPhi_Error}, \ref{RT0_DPhi_Error}, \ref{BDM1_DPhi_Error}, \ref{RT1_DPhi_Error},
\ref{BDM2_DPhi_Error} show the errors for $|\nabla\phi|_2$ using a $128^2$ mesh for solving the boundary of Fig.
\ref{INTERFACE1}.

In Table \ref{INTERFACE_ERROR}, we show the maximum gradient errors on the boundary by different methods. From this table, we know that the order of accuracies of maximum gradient errors on the boundary for the five methods are comparable with the $L_2$ norm on the whole domain.
\begin{table}
\caption{Maximum gradient errors on the boundary by different methods}
\label{INTERFACE_ERROR}
\begin{tabular}{|c||r|r|r|r|r|}
\hline
 Size           &  EBM          &  RT0          &  BDM1         &RT1            & BDM2\\
\hline
$64\times64$    &2.459720e-03   & 9.367834e-03  &7.927028e-04   & 3.102073e-04  & 4.988912e-05      \\
$128\times128$  &6.567893e-04   & 6.797467e-03  &1.274594e-04   & 4.007023e-05  & 1.038527e-06      \\
$256\times256$  &1.755489e-04   & 3.626596e-03  &2.486447e-05   & 1.321007e-05  & 1.033665e-07     \\
$512\times512$  &4.614643e-05   & 1.754357e-03  &6.351849e-06   & 2.751578e-06  & 2.310828e-08     \\
\hline
\end{tabular}
\end{table}

\subsection{Automatic Mesh Refinement vs. Uniform Grid}

\begin{figure}[htp]
\centering
\includegraphics[scale=0.3]{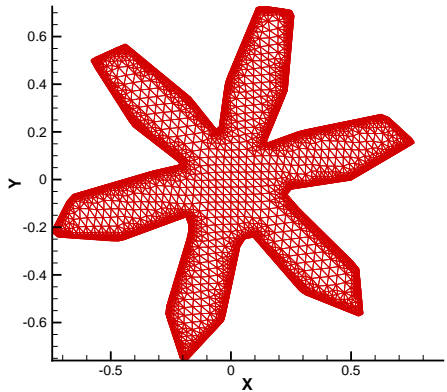}
\caption{Mesh when minimum level is 6 and maximum level is 9}
\label{Mesh_6x9}
\end{figure}

\begin{figure}[htp]
\centering
\includegraphics[scale=0.3]{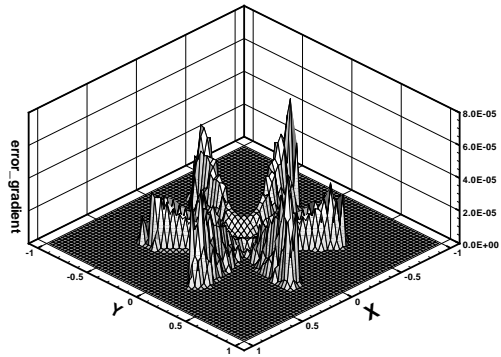}
\caption{Flux error when minimum level is 6 and maximum level is 9}
\label{RT1_DPhi_Error_6x9}
\end{figure}

In this section, we compare the convergence rates when the mesh is refined around the boundary. The boundary is the same as the boundary of Fig. \ref{INTERFACE1} for our second test. We only compute the results using RT1. Figure \ref{Mesh_6x9} shows the mesh when the quadtree's minimum level is 6 and maximum level is 9. And Fig. \ref{RT1_DPhi_Error_6x9} gives the flux error plot. From Table \ref{TABLE_CONV_TIMING_AMR} we see that the refinement does not give a more accurate solution on the whole computational domain. In fact, this is reasonable. From Fig. \ref{RT1_DPhi_Error_6x9}, we know that the maximum errors are located in the interior where the mesh is coarsest, where the mesh is not refined. Our refinement is only around the boundary and in this way the errors near the boundary are reduced. The effect of boundary flux error reduction by mesh refinement diminish gradually.

\begin{table}
\caption{Convergence and timing results for automatic mesh refinement 
(minimum level = 6) using RT1 for the domain of Fig. \ref{INTERFACE1} }
\label{TABLE_CONV_TIMING_AMR}
\begin{tabular}{|c||l|l|l|l|l|}
\hline
maximum         &\multicolumn{4}{c|}{RT1} &{unknown}\\
\cline{2-5}
level          &  error         & max boundary error   &  time    &iterations  &number\\
\cline{2-6}
6              &2.203159e-05    &3.102045e-04       &1.461276  & 151        & 6354\\
7              &1.766854e-05    &4.128689e-05       &2.533742  & 224        & 11200\\
8              &1.539364e-05    &8.678919e-06       &7.155402  & 330        & 22758\\
9              &1.462538e-05    &3.401366e-06       &20.458581 & 467        & 47214\\
10             &1.387916e-05    &8.151048e-07       &68.288392 & 789        & 96076\\
\hline
\end{tabular}
\end{table}


\section{Conclusion}

In this paper, we have used the embedded boundary method and the mixed finite element method to solve the elliptic boundary value problem in 2D. We compared the convergence and timing results. 

Since the embedded boundary method uses a structured cartesian grid, it is easier to implement. It is much harder to write the mesh generation program. But after the mesh is given, the discretization is simpler for the mixed finite element method. And it is easier to use the mixed finite element for the elliptic interface problem since the interface is in fact an internal boundary. However, the EBM method must be modified to solve an elliptic interface problem. To save computational resources when solving large problems, we could use EBM with automatic mesh refinement, which is one important part of our mesh generation method. 

The EBM has the advantage of fewer unknowns with the same mesh size compared with the MFM. There are two reasons for this. One reason is that the EBM uses a structured grid and the finite volume/central finite difference has super convergence in the mesh. The MFM uses an unstructured grid, and to achieve the same order of accuracy, a higher order basis function space is needed, which means more unknowns. The other reason is that the unknowns for EBM are cell centered and those for the MFM are edge centered. Since the approximate ratio of the vertices to faces to edges is 1:2:3 for a simple large triangle mesh, we know the ratio of the unknowns for the EBM, RT0, BDM1, RT1, BDM2 is approximately 1:3:6:6:9. Thus the EBM problem is smaller, which explains why it is much more faster. However, for a given accuracy, the fastest method is BDM2 which is 3rd order accurate in flux, and then EBM/RT1/BDM1/RT0.

\bibliographystyle{plain}
\bibliography{xBib}

\begin{thebibliography}{10}

\bibitem{MckenneyGreenGardMayo1995}
A.~Mayo A.~Mckenney, L.~Greengard.
\newblock A fast poisson solver for complex geometries.
\newblock {\em J. Comput. Phys.}, 118:348每355, 1995.

\bibitem{Chavent1999}
P.~Ackerer A.~Younes, R.~Mose and G.~Chavent.
\newblock A new formulation of the mixed finite element method for solving
  elliptic and parabolic pde with triangular elements.
\newblock {\em J. Comput. Phys.}, 149(1):148--167, 1999.

\bibitem{Arnold1990}
Douglas~N. Arnold.
\newblock Mixed finite element methods for elliptic problems.
\newblock {\em Comput. Methods Appl. Mech. Engrg.}, 82:281--300, 1990.

\bibitem{BrezziFortin1991}
Franco Brezzi and Michel Fortin.
\newblock {\em Mixed and Hybrid Finite Element Methods}.
\newblock Springer Series In Computational Mathematics 15, 1991.

\bibitem{Chavent1991}
G.~Chavent and J.E. Roberts.
\newblock A unified physical presentation of mixed, mixed-hybrid finite
  elements and standard finite difference approximations for the determination
  of velocities in waterflow problems.
\newblock {\em Adv. Water Resources}, 14(6):329--348, 1991.

\bibitem{JohansenColella98}
P.~Colella H.~Johansen.
\newblock A cartesian grid embedding boundary method for poisson's equation on
  irregular domains.
\newblock {\em J. Comput. Phys.}, 147:60每85, 1998.

\bibitem{LeVequeLi1994}
R.J. LeVeque and Z.L. Li.
\newblock The immersed interface method for elliptic equations with
  discontinuous coefficients and singular sources.
\newblock {\em SIAM J. Numer. Anal.}, 31:1019每1044, 1994.

\bibitem{Marcum2001}
D.~L. Marcum.
\newblock Efficient generation of high-quality unstructured surface and volume
  grids.
\newblock {\em Engineering with Computers}, 17(3):211--233, 2001.

\bibitem{YerryShephard1983}
M.S.Shephard M.A.Yerry.
\newblock A modified quadtree approach to finite element generation.
\newblock {\em IEEE Comput. Graph. Appl.}, 3(1):39--46, 1983.

\bibitem{Mayo1984}
A.~Mayo.
\newblock The fast solution of poisson's and the biharmonic equations on
  irregular regions.
\newblock {\em SIAM J. Numer. Anal.}, 21:285每299, 1984.

\bibitem{Berg1998}
M.~Overmars M.D.~Berg, M.V.~Kreveld and O.~Schwarzkopf.
\newblock {\em Computational Geometry: Algorithms and Applications}.
\newblock Springer, second edition edition, 1998.

\bibitem{Peskin2002}
Charles~S. Perskin.
\newblock The immersed boundary method.
\newblock {\em Acta Numerica}, 11:479--517, 2002.

\bibitem{Samet1990}
Hanan Samet.
\newblock {\em The design and analysis of spatial data structures}.
\newblock Addison-Wesley Series In Computer Science, 1990.

\bibitem{Shewchuk1997}
Jonathan~Richard Shewchuk.
\newblock {\em Delaunay Refinement Mesh Generation}.
\newblock PhD thesis, Computer Science Department, Carnegie Mellon University,
  May 1997.

\bibitem{Wang2006}
Shuqiang Wang.
\newblock {\em A Mixed Finite Element Method for Elliptic Interface/Boundary
  Value Problem}.
\newblock PhD thesis, SUNY at Stony Brook, August 2006.

\bibitem{LorensenCline1987}
H.E.~Cline W.E.~Lorensen.
\newblock Marching cubes: a high resolution 3d surface construction algorithm.
\newblock {\em Computer Graphics}, 21(4):163--169, 1987.

\bibitem{LiuFedkiwKang2000}
M.~Kang X.D.~Liu, R.P.~Fedkiw.
\newblock A boundary condition capturing method for poisson's equation on
  irregular domains.
\newblock {\em J. Comput. Phys.}, 160, 2000.

\end{thebibliography}

\end{document}